\documentclass[]{siamart171218}

\usepackage{algorithmic}
\usepackage{amsmath}
\usepackage{amssymb}
\usepackage{graphicx}
\usepackage{listings}

\newcommand{\X}{{\bf x}}

\newcommand{\DX}{{\Delta \bf x}}
\newcommand{\DdX}{{\Delta \dot{\bf x}}}

\newcommand{\DF}{{\Delta F}}

\newcommand{\B}{{\bf b}}
\newcommand{\C}{{\bf c}}

\newcommand{\DZ}{{\Delta \bf z}}
\renewcommand{\P}{{\bf p}}

\newcommand{\Y}{{\bf y}}
\newcommand{\A}{{\bf a}}

\newcommand{\ass}{\small \mathrel{{:}{=}}}
\newcommand{\R}{\mathbb{R}}
\newcommand{\Z}{{\bf z}}
\renewcommand{\r}{{\bf r}}

\title{Numerical Stability 
of Tangents and Adjoints of Implicit Functions}
\author{Uwe Naumann\thanks{Informatik 12: Software and Tools for Computational Engineering, RWTH Aachen University, Germany. \email{naumann@stce.rwth-aachen.de}}}

\begin{document}

\maketitle

%\begin{keywords} algorithmic differentiation, tangent, adjoint \end{keywords}
\begin{abstract}
We investigate errors in tangents and adjoints of implicit functions
resulting from errors in the primal solution due to
approximations computed by a numerical solver. 

Adjoints of systems of linear equations turn out to be unconditionally 
numerically stable. Tangents of systems of linear equations can become 
instable as well as both tangents and adjoints of systems of nonlinear 
equations, which extends to optima of convex unconstrained objectives. 
Sufficient conditions for numerical stability are derived.
\end{abstract}

\section{Introduction}

We consider twice differentiable implicit functions
\begin{equation} \label{eqn:F}
	F : \R^m \rightarrow \R^n: \P \mapsto \X=F(\P) 
\end{equation} 
defined by the roots of residuals
\begin{equation} \label{eqn:R}
	R : \R^n \times \R^m \rightarrow \R^n: (\X,\P) \mapsto R(\X,\P) \; .
\end{equation} 
$R$ is referred to as the primal residual as opposed to tangent and adjoint 
residuals to be considered later.
Primal roots of the residual satisfying 
\begin{equation} \label{eqn:R0}
R(\X,\P)=0 
\end{equation} 
are assumed to be approximated by numerical solvers
$$
S : \R^m \rightarrow \R^n: \P \mapsto \X+\Delta \X=S(\P) 
$$
with an absolute error $\Delta \X$ yielding a 
relative error $\delta \X$ of norm
$$
\|\delta \X\|= \frac{\|\Delta \X\|}{\|\X\|} = \frac{\|S(\P)-F(\P)\|}{\|F(\P)\|}\; .
$$
We investigate (relative) errors in corresponding tangents
\begin{equation} \label{eqn:tan}
\dot{\X} = \dot{F}(\X,\dot{\P}) \equiv 
\frac{d F}{d \P} \cdot \dot{\P}
\end{equation} 
and adjoints
\begin{equation} \label{eqn:adj}
\bar{\P} = \bar{F}(\X,\bar{\X}) \equiv 
\frac{d F}{d \P}^T \cdot \bar{\X}
\end{equation} 
due to $\Delta \X.$ 
{\em Algorithmic} tangents and adjoints result from the application of 
algorithmic differentiation (AD) \cite{Griewank2008EDP,Naumann2012TAo} to the
solver $S.$ 
{\em Symbolic} tangents and adjoints can be 
derived at the solution of Equation~(\ref{eqn:R0}) in terms of tangents and 
adjoints of the residual \cite{Giles2008CMD,Naumann2015NLS}. 
AD of the solver can thus be avoided which typically
results in a considerably lower computational complexity.

\section{Prerequisites}
We perform standard first-order error analysis. 
For a given absolute error $\Delta \P$ in the input of a function $F$ 
the absolute error in the result is estimated as
\begin{equation}  \label{eqn:foe}
\Delta \X \approx \frac{d F}{d \P} \cdot \Delta \P  \; .
\end{equation}
Equation~(\ref{eqn:F}) is differentiated with respect to $\P$ in 
the direction of the absolute error $\Delta \P.$ From the Taylor series
expansion  of 
$$
\X + \Delta \X = \X +  \frac{d F}{d \P} \cdot \Delta \P + O(\|\Delta \P\|^2) 
$$
it follows that
negligence of the remainder within a neighborhood of $\X$ containing 
$\Delta \X$ is reasonable for $\|\Delta \P\| \rightarrow 0$ and assuming
convergence of the Taylor series to the correct function value. 
%Equation~(\ref{eqn:foe}) follows immediately.
For linear $F$ we get
$\Delta \X =\frac{d F}{d \P} \cdot \Delta \P$ due to the vanishing remainder.

%First-order tangent error analysis yields 
%\begin{equation} \label{eqn:tanerr}
%\|\delta \dot \X\| =
%	\frac{\|\Delta \dot \X\|}{\|\dot \X\|} = 
%	\frac{\|\Delta \dot F(\X,\dot{\P})\|}{\|\dot F(\X,\dot \P)\|}  
%	\approx
%\frac{\|\dot F_\X(\X,\dot \P)\|}{\|\dot F(\X,\dot \P)\|} \cdot \|\Delta \X\| \;.
%\end{equation} 
%for any submultiplicative norm $\|.\|.$ 
%Similarly, the adjoint error is estimated as
%\begin{equation} \label{eqn:adjerr}
%\|\delta \bar \P\|=\frac{\|\Delta \bar \P\|}{\|\bar \P\|}
%	=\frac{\|\Delta \bar F(\X,\bar \X)\|}{\|\bar F(\X,\bar \X)\|}  
%	\approx
%\frac{\|\bar F_\X(\X,\bar \X)\|}{\|\bar F(\X,\bar \X)\|} \cdot \|\Delta \X\| \;.
%\end{equation} 
%

Tangents and adjoints of Equation~(\ref{eqn:F}) can be expressed as matrix 
equations over derivatives of the residual. The fundamental operations 
involved are scalar multiplications and additions, outer vector 
products, matrix-vector products and solutions of systems of linear equations. 
In this section we recall the corresponding well-known first-order error 
estimates. Norms of those estimates are considered.

\subsection{Scalar Multiplication}

Differentiation of 
$y=x_1 \cdot x_2$ in direction 
$$\Delta \X = \begin{pmatrix} x_1 \\ x_2 \end{pmatrix} \neq 0$$ 
yields the absolute error
$|\Delta y|=|\Delta x_1 \cdot x_2+x_1 \cdot \Delta x_2|.$ 
The corresponding relative error is equal to
\begin{align*}
|\delta y| \equiv
\frac{|\Delta y|}{|y|} &
	= \frac{|x_2 \cdot \Delta x_1 +x_1 \cdot \Delta x_2|}{|x_1 \cdot x_2|} 
	= \left | \frac{x_2 \cdot \Delta x_1 +x_1 \cdot \Delta x_2}{x_1 \cdot x_2} \right |\\
	&= \left | \frac{\Delta x_1}{x_1} + \frac{\Delta x_2}{x_2} \right |
	= \left | \frac{\Delta x_1}{x_1} \right | + \left |\frac{\Delta x_2}{x_2} \right |
	= \frac{|\Delta x_1|}{|x_1|} + \frac{|\Delta x_2|}{|x_2|} 
	= |\delta x_1|+|\delta x_2| \; .
\end{align*}
Scalar multiplication turns out to be numerically stable. The relative error
in the result is of the order of the maximum relative error in the arguments.
A similar result holds for scalar division.
It generalizes naturally to element-wise multiplication and division of 
vectors, matrices, and higher-order tensors as well as to the outer product
of two vectors.

\subsection{Scalar Addition} 
Differentiation of 
$y=x_1+x_2$ in direction $\Delta \X \neq 0$ yields the absolute error
$|\Delta y|=|\Delta x_1+\Delta x_2|.$ The corresponding relative error is equal
to
$$
|\delta y| \equiv
\frac{|\Delta y|}{|y|}= 
\frac{|\Delta x_1 +\Delta x_2|}{|x_1+x_2|} \; .
$$
Scalar addition turns out to be numerically unstable due to 
$|\delta y| \rightarrow \infty$ for $\Delta x_1\neq -\Delta x_2$ and $x_1 \rightarrow -x_2.$ 
The relative error in the result can become unbounded for bounded relative
errors of the arguments.
A similar result holds for scalar subtraction.

Numerical instability of scalar addition prevents unconditional 
numerical stability of inner vector products as well as matrix-vector/matrix 
products and solutions of systems of linear equations. Sufficient conditions for
numerical stability need to be formulated.

\subsection{Matrix-Vector Product} \label{sec:Mxv}
Differentiation of the matrix-vector product
$\X=A \cdot \B$ for
$A \in \R^{n \times n}$ and $\X,\B \in \R^n$ 
in the direction of non-vanishing absolute errors
$\Delta A \in \R^{n \times n}$ and $\Delta \X, \Delta \B \in \R^n$ 
yields 
\begin{align*}
	\Delta \X &= \Delta A \cdot \B + A \cdot \Delta \B  
\end{align*}
and hence the first-order error estimate
\begin{align*}
	\frac{\|\Delta \X\|}{\|\X\|} &= \frac{\|\Delta A \cdot \B + A \cdot \Delta \B \|}{\|A \cdot \B\|} \\
	&\leq \frac{\|\Delta A\| \cdot \|\B\| + \|A\| \cdot \|\Delta \B \|}{\|A \cdot \B\|} 
	= \|A^{-1}\| \cdot \frac{\|\Delta A\| \cdot \|\B\| + \|A\| \cdot \|\Delta \B \|}{\|A^{-1}\| \cdot \|A \cdot \B\|} \\
	&\leq \|A^{-1}\| \cdot \frac{\|\Delta A\| \cdot \|\B\| + \|A\| \cdot \|\Delta \B \|}{\|A^{-1} \cdot A \cdot \B\|} 
	= \|A^{-1}\| \cdot \frac{\|\Delta A\| \cdot \|\B\| + \|A\| \cdot \|\Delta \B \|}{\|\B\|} \\
	&= \|A^{-1}\| \cdot \|A\| \cdot \left ( \frac{\|\Delta A\|}{\|A\|} + \frac{\|\Delta \B \|}{\|\B\|} \right )
	= \kappa(A) \cdot \left (\frac{\|\Delta A\|}{\|A\|} + \frac{\|\Delta \B \|}{\|\B\|} \right )\; .
\end{align*}
A low condition number 
$\kappa(A) \equiv \|A^{-1}\| \cdot \|A\|$ of $A$ is sufficient for numerical 
stability. We write
\begin{equation} \label{err:Mxv}
\|\delta \X\| \approx \kappa(A) \cdot (\|\delta A\| + \|\delta \B \|) 
\end{equation}
to indicate that depending on the magnitude of $\kappa(A)$
the relative error of a matrix-vector product can suffer from 
a potentially dramatic amplification of the relative errors in the 
arguments. 

\subsection{Systems of Linear Equations} \label{sec:Axb}

Differentiation of the system of linear equations
$A \cdot \X=\B$ for
$A \in \R^{n \times n}$ and $\X,\B \in \R^n$ 
in the direction of non-vanishing absolute errors
$\Delta A \in \R^{n \times n}$ and $\Delta \X, \Delta \B \in \R^n$ 
yields 
\begin{align*}
\Delta \X= A^{-1} \cdot (\Delta \B - \Delta A \cdot \X)
\end{align*}
and hence the first-order error estimate
\begin{align*}
	\frac{\|\Delta \X\|}{\|\X\|} &= \frac{\|A^{-1} \cdot (\Delta \B - \Delta A \cdot \X) \|}{\|\X\|} \\
	&\leq \frac{\|A^{-1} \cdot \Delta \B\|}{\|\X\|}+\frac{\| A^{-1} \cdot \Delta A \cdot \X \|}{\|\X\|} \\
	&\leq \frac{\|A^{-1}\| \cdot \|\Delta \B\|}{\|\X\|}+\frac{\| A^{-1} \cdot \Delta A \cdot \X \|}{\|\X\|} 
	= \frac{\|A\| \cdot \|A^{-1}\| \cdot \|\Delta \B\|}{\|A\| \cdot \|\X\|}+\frac{\| A^{-1} \cdot \Delta A \cdot \X \|}{\|\X\|} \displaybreak[1] \\
	&\leq \kappa(A) \cdot \frac{\|\Delta \B\|}{\|A \cdot \X\|}+\frac{\| A^{-1} \cdot \Delta A \cdot \X \|}{\|\X\|} 
	= \kappa(A) \cdot \frac{\|\Delta \B\|}{\|\B\|}+\frac{\| A^{-1} \cdot \Delta A \cdot \X \|}{\|\X\|} \displaybreak[1] \\
	&\leq \kappa(A) \cdot \frac{\|\Delta \B\|}{\|\B\|}+\frac{\| A^{-1}\| \cdot \|\Delta A\| \cdot \|\X \|}{\|\X\|} =
	\kappa(A) \cdot \frac{\|\Delta \B\|}{\|\B\|}+\| A^{-1}\| \cdot \|\Delta A\| \\
	&= \kappa(A) \cdot \frac{\|\Delta \B\|}{\|\B\|}+\frac{\|A\| \cdot \| A^{-1}\| \cdot \|\Delta A\|}{\|A\|} 
	= \kappa(A) \cdot \frac{\|\Delta \B\|}{\|\B\|}+\kappa(A) \cdot \frac{\|\Delta A\|}{\|A\|} \\
	&= \kappa(A)\cdot \left ( \frac{\|\Delta \B\|}{\|\B\|}+ \frac{\|\Delta A\|}{\|A\|} \right ) \; .
\end{align*}
As in Section~\ref{sec:Mxv} we get
\begin{equation} \label{err:Axb}
\|\delta \X\| \approx \kappa(A) \cdot (\|\delta A\| + \|\delta \B \|) \; . 
\end{equation}
Again, a low condition number of $A$ is sufficient for numerical 
stability. 

\section{Errors in Tangents and Adjoints of Implicit Functions}

Differentiation of Equation~(\ref{eqn:R0}) with respect to $\P$ yields
\begin{equation} \label{eqn:dR}
	 \frac{\partial R}{\partial \X} \cdot \frac{d \X}{d \P} + \frac{\partial R}{\partial \P} = R_\X \cdot \frac{d \X}{d \P} + R_\P=0 \; ,
\end{equation} 
where $\partial$ denotes partial differentiation. Multiplication with $\dot \P$ from the right yields the tangent residual
\begin{equation} \label{eqn:tanR}
	R_\X \cdot \frac{d \X}{d \P} \cdot \dot \P + R_\P \cdot \dot \P=
	R_\X \cdot \dot \X + R_\P \cdot \dot \P =0 \; .
\end{equation} 
The tangent $\dot \X$ can be computed as the solution of the system of 
linear equations 
$$
	R_\X \cdot \dot \X =- R_\P \cdot \dot \P \; .
$$
The right-hand side is obtained by a single evaluation of the tangent residual.
Tangents in the directions of the Cartesian basis of $\R^n$ yields $R_\X.$ 
Potential sparsity can and should be exploited \cite{Gebremedhin2005WCI}.
An error $\Delta \X$ in the primal solution yields a 
corresponding error in the tangent for $R_\X=R_\X(\X)$ and/or $R_\P=R_\P(\X).$

From Equation~(\ref{eqn:dR}) it follows that for regular $R_\X$
$$
\frac{d \X}{d \P}= - R_\X^{-1} \cdot R_\P \; .
$$
Transposition of the latter followed by multiplication with $\bar \X$ from the right yields
\begin{equation} \label{eqn:adjR}
	\bar \P = \frac{d \X}{d \P}^T \cdot \bar \X = - R_\P^T \cdot R_\X^{-T} \cdot \bar \X \; .
\end{equation} 
The adjoint $\bar \P$ can be computed as the solution of the system of 
linear equations $$ R_\X^T \cdot \Z = -\bar \X$$ followed by the evaluation of
the adjoint residual yielding $$\bar \P=R_\P^T \cdot \Z$$. Again, an error 
$\Delta \X$ in the primal solution yields a 
corresponding error in the adjoint.

\subsection{Systems of Linear Equations}

The tangent of 
the solution of the primal system of linear equations
\begin{equation} \label{primal_ls}
A \cdot \X = \B
\end{equation}
is defined as $\dot{\X} = \dot{\X}_A + \dot{\X}_\B,$ where  
\begin{equation} \label{eqn:tls1}
A \cdot \dot{\X}_\B = \dot{\B}
\end{equation}
and
\begin{equation} \label{eqn:tls2}
A \cdot \dot{\X}_A = -\dot{A} \cdot \X 
\end{equation}
\cite{Giles2008CMD}.
An error $\Delta \X$ in the primal solution which, for example, might result from the use of 
an indirect solver yields an erroneous tangent 
$$
\dot{\X} + \Delta \dot{\X}= (\dot{\X}_A + \Delta \dot{\X}_A) + (\dot{\X}_\B + \Delta \dot{\X}_\B) \; .
$$
Application of Equation~(\ref{err:Axb}) to Equation~\eqref{eqn:tls1} yields
$$
	\|\delta \dot{\X}_\B\| \approx \kappa(A)\cdot ( \|\delta A\| + \|\delta \dot{\B}\|) \; .
$$
Independence of $\dot{\X}_\B$ from $\X$ (and hence from $\Delta \X$) 
implies $\delta \dot{\X}_\B=0$ for error-free $A$ and $\dot{\B},$ that is
$
\Delta \dot{\X} = \Delta \dot{\X}_A$, respectively
$
\delta \dot{\X} = \delta \dot{\X}_A.
$
Let $\C=-\dot A \cdot \X.$ With Equation~(\ref{err:Mxv}) it follows that
$$
\|\delta \C\| \approx \kappa(\dot A) \cdot \|\delta \X\| 
$$
as $\delta \dot A=0$
Moreover, application of Equation~(\ref{err:Axb}) to $A \cdot \dot \X_A = \C$ yields
$$
\|\delta \dot \X_A\| \approx \kappa(A) \cdot \|\delta \C\| \; .
$$
Consequently, 
\begin{equation} \label{err:tls}
\|\delta \dot \X_A\| \approx \kappa(A) \cdot \kappa(\dot A) \cdot \|\delta \X\| \; .
\end{equation}
Low condition numbers of both $A$ and $\dot A$ ensure
numerical stability of tangent systems of linear equations.

The adjoint of the primal linear system in Equation~(\ref{primal_ls})
is defined as 
\begin{equation} \label{eqn:als1} 
A^T \cdot \bar{\B} = \bar{\X}
\end{equation}
and 
\begin{equation} \label{eqn:als2} 
\bar{A} = -\bar{\B} \cdot \X^T  
\end{equation}
\cite{Giles2008CMD}.
Application of Equation~(\ref{err:Axb}) to Equation~(\ref{eqn:als1}) yields
$$
	\delta \bar{\B} \approx \kappa(A)\cdot ( \delta A + \delta \bar{\X}) \; .
$$
Independence of $\bar \B$ from $\X$ (and hence from $\Delta \X$) 
implies $\delta \bar{\B}=0$ for error-free $A$ and $\bar{\X}.$ 
The outer product
$\bar{A} = -\bar{\B} \cdot \X^T$ is numerically stable as
scalar multiplication is. Consequently, adjoint systems of linear equations
are numerically stable.

\subsection{Systems of Nonlinear Equations}

Differentiation of
Equation~(\ref{eqn:tanR}) in the direction of absolute errors 
$\Delta R_\X \in \R^{n \times n},$ $\Delta \dot{\X} \in \R^n,$ $\Delta R_\P \in \R^{n \times m}$ and $\Delta \dot{\P} \in \R^m$
yields
$$
\Delta R_\X \cdot \dot{\X} + R_\X \cdot \Delta \dot{\X} +\Delta R_\P \cdot \dot{\P}~[+\underset{=0}{\underbrace{R_\P \cdot \Delta \dot{\P}}}]=0 
$$
as $\Delta \dot{\P}=0$ and hence
\begin{align*}
	\Delta \dot{\X}& = R^{-1}_\X \cdot (\Delta R_\X \cdot \dot{\X} + \Delta R_\P \cdot \dot{\P} )  \; .
\end{align*}
First-order estimates for 
$$\Delta R_\X \cdot \dot{\X}=[\Delta R_\X \cdot \dot{\X}]_i \approx [R_{\X,\X}]_{i,j,k} \cdot  [\dot{\X}]_j \cdot [\Delta \X]_k \equiv \Delta \dot R_\X \cdot \Delta \X$$ 
and 
$$\Delta R_\P \cdot \dot{\P}=[\Delta R_\P \cdot \dot{\P}]_i \approx [R_{\P,\X}]_{i,j,k} \cdot [\dot{\P}]_j \cdot [\Delta \X]_k
\equiv \Delta \dot R_\P \cdot  \Delta \X$$ 
in index notation
(summation over the shared index) yield
\begin{align*}
	\Delta \dot{\X}
	& \approx R^{-1}_\X \cdot (\Delta \dot R_\X + \Delta \dot R_\P ) \cdot \Delta \X 
\end{align*}
and hence, with Equation~(\ref{err:Mxv}),
\begin{equation} \label{err:tnls}
	\|\delta \dot{\X}\| \approx
	\kappa(R_\X) \cdot \kappa(\Delta \dot R_\X + \Delta \dot R_\P) 
	\cdot \|\delta \X\| \; .
\end{equation}
Low condition numbers 
of the respective first and second derivatives of the residual ensure
numerical stability of tangent systems of nonlinear equations.
Both $\Delta \dot R_\X$ and $\Delta \dot R_\P$ can be computed by 
algorithmic differentiation (AD) \cite{Griewank2008EDP,Naumann2012TAo}.

Application of Equation~(\ref{err:Axb}) to the system of linear equations
$$
R_\X^T \cdot \Z = - \bar \X
$$
for $\Delta \bar \X=0$
yields
$$
\Delta \Z =R_\X^{-T} \cdot \Delta R_\X^T \cdot \Z 
$$
and hence
$$
\|\delta \Z \| \approx \kappa(R_\X) 
	\cdot \kappa(\Delta \bar R_\X) 
	\cdot \|\delta \X\| \; ,
$$
where 
$$[\Delta R^T_\X \cdot \Z]_j \approx [R_{\X,\X}]_{i,j,k} \cdot  [\Z]_i \cdot [\Delta \X]_k \equiv \Delta \bar R_\X \cdot \Delta \X \; .$$

Differentiation of $\bar \P=R_\P^T \cdot \Z$ in the direction of
the non-vanishing absolute errors $\Delta R_\P^T \in \R^{m \times n}$ and $\Delta \Z \in \R^n$
yields
\begin{align*}
	\Delta \bar \P &= \Delta R_\P^T \cdot \Z + R_\P^T \cdot \Delta \Z
= \Delta R_\P^T \cdot \Z + R_\P^T \cdot 
R_\X^{-T} \cdot \Delta R_\X^T \cdot \Z 
 \end{align*}
and hence
\begin{equation} \label{err:anls}
	\|\delta \bar{\P}\| \approx
	\left ( \kappa(\Delta \bar R_\P) + 
	\kappa(R_\P) \cdot 
	\kappa(R_\X) \cdot 
	\kappa(\Delta \bar R_\X) \right )
	\cdot \|\delta \X\| \; ,
\end{equation}
where 
$$[\Delta R^T_\P \cdot \Z]_j \approx [R_{\P,\X}]_{i,j,k} \cdot  [\Z]_i \cdot [\Delta \X]_k \equiv \Delta \bar R_\P \cdot \Delta \X \; .$$
Low condition numbers 
of the respective first and second derivatives of the residual ensure
numerical stability of adjoint systems of nonlinear equations.
Both $\Delta \bar R_\X$ and $\Delta \bar R_\P$ can be computed by AD.

\subsection{Convex Unconstrained Objectives}
The first-order optimality condition 
for a parameterized convex unconstrained objective
$$
f : \R^n \times \R^m \rightarrow \R: (\X,\P) \mapsto y=f(\X,\P) 
$$
yields the residual
$f_\X(\X,\P)=0.$ 
Consequently, assuming $f$ to be three times differentiable,
\begin{equation} \label{err:tnlp}
	\|\delta \dot{\X}\| \approx
	\kappa(f_{\X,\X}) \cdot \kappa(\Delta \dot f_{\X,\X} + \Delta \dot f_{\X,\P})
	\cdot \|\delta \X\|  \; ,
\end{equation} 
where
$$\Delta f_{\X,\X} \cdot \dot{\X}=[\Delta f_{\X,\X} \cdot \dot{\X}]_i \approx [f_{\X,\X,\X}]_{i,j,k} \cdot  [\dot{\X}]_j \cdot [\Delta \X]_k \equiv \Delta \dot f_{\X,\X} \cdot \Delta \X$$ 
and 
$$\Delta f_{\X,\P} \cdot \dot{\P}=[\Delta f_{\X,\P} \cdot \dot{\P}]_i \approx [f_{\X,\P,\X}]_{i,j,k} \cdot [\dot{\P}]_j \cdot [\Delta \X]_k
\equiv \Delta \dot f_{\X,\P} \cdot \Delta \X \; .$$ 
Similarly,
\begin{equation} \label{err:anlp}
\|\delta \bar{\P}\| \approx
	\left (\kappa(\Delta \bar f_{\X,\P}) + 
	\kappa(f_{\X,\P}) \cdot 
	\kappa(f_{\X,\X}) \cdot 
	\kappa(\Delta \bar f_{\X,\X}) \right )
	\cdot \|\delta \X\| \; ,
\end{equation} 
where
$$[\Delta f_{\X,\X}^T \cdot \Z]_j = [\Delta f_{\X,\X} \cdot \Z]_j \approx [f_{\X,\X,\X}]_{i,j,k} \cdot  [\Z]_i \cdot [\Delta \X]_k \equiv \Delta \bar f_{\X,\X} \cdot \Delta \X $$
and 
$$[\Delta f_{\X,\P}^T \cdot \Z]_j \approx [f_{\X,\P,\X}]_{i,j,k} \cdot  [\Z]_i \cdot [\Delta \X]_k \equiv \Delta \bar f_{\X,\P} \cdot \Delta \X \; .$$
Low 
condition numbers of the respective second and third derivatives of the
objective ensure
numerical stability of tangent and adjoint optima of convex unconstrained 
objectives.
Both $\Delta \dot f_{\X,\X}$ and $\Delta \dot f_{\X,\P}$ as well as
$\Delta \bar f_{\X,\X}$ and $\Delta \bar f_{\X,\P}$ can be computed by AD.

\section{Conclusion}

Adjoint systems of linear equations are numerically stable with respect to 
errors in the primal solution. However, numerical stability of tangents and 
adjoints of implicit functions cannot be 
guaranteed in general. Sufficient conditions in terms of derivatives of the 
residual are given by Equations~(\ref{err:tls}), (\ref{err:tnls}), 
(\ref{err:anls}), (\ref{err:tnlp}) and (\ref{err:anlp}). AD
can be used to compute 
these derivatives. Corresponding symbolic tangents and adjoints 
should be 
augmented with optional estimation of conditions of the relevant derivatives.

\end{document}